\documentclass[11pt]{amsart}
\usepackage{amssymb,amsmath,amsthm}
\usepackage[margin=1.1in]{geometry}
\usepackage{enumitem}
\usepackage{tikz}
\usepackage{hyperref}

\newtheorem{theorem}{Theorem}[section]
\newtheorem{lemma}[theorem]{Lemma}
\newtheorem{proposition}[theorem]{Proposition}
\newtheorem{corollary}[theorem]{Corollary}
\newtheorem{thmC}{Main Theorem}

\newtheorem{thmB}{Main Theorem}

\newtheorem{thmA}{Main Theorem}

\theoremstyle{definition}
\newtheorem{definition}[theorem]{Definition}
\newtheorem{remark}[theorem]{Remark}

\newcommand{\lk}{\operatorname{lk}}
\newcommand{\Z}{\mathbb Z}
\newcommand{\F}{\mathbb F}
\newlength{\faclen}
\newenvironment{faclist}[1]{\settowidth{\faclen}{\small\ttfamily #1}\par\smallskip\noindent\hfil\begin{minipage}{\faclen}\raggedright\small\ttfamily}{\end{minipage}\hfil\par\smallskip}

\title{A constructible non-shellable $3$-sphere}
\author{Bruno Benedetti}
\address[BB]{Department of Mathematics, University of Miami, Coral Gables, FL, USA}
\email{bruno@math.miami.edu}
\author{Giuliamaria Menara}
\address[GM]{CNRS, Institut de Math\'ematiques de Jussieu-Paris Rive Gauche, Paris, France}
\email{giuliamaria.menara@imj-prg.fr}
\date{\today}

\begin{document}
\begin{abstract}
We extend the obstructions to shellability, constructibility and vertex-decomposability from knots to links. We also refine the existing obstructions to shellability, by following how a link grows along a shelling, in the spirit of persistent homology. This allows us to give the first example of a constructible, unshellable $3$-sphere, thereby settling a long-standing problem in combinatorial topology.
\end{abstract}
\maketitle

\section{Introduction}

For triangulated $3$-spheres, vertex-decomposable implies shellable implies constructible. Whether the second implication is strict was called a ``big open question'' in~\cite{HS}, and later relaunched in~\cite{BZ,BL,GG,Ven}, among others. The problem is that for spheres, the only known obstructions to these desirable combinatorial properties come from knots in the $1$-skeleton: a knot with $e$ edges and bridge number $b$ forces $e\ge2b$ in a shellable or constructible sphere, and $e\ge3b$ in a vertex-decomposable one~\cite{EH,HS}. Both bounds are tight. Thus the difficulty in finding a constructible unshellable sphere is that the only known obstruction to shellability prevents constructibility as well.

To bypass this problem, in this paper we study links, i.e.\ disjoint unions of knots. We extend all obstructions from knots to links, and we find out that for shellability one can be more precise.

\begin{thmA}[Proposition~\ref{prop:linkconstr} and Theorems~\ref{lem:linkbound}, \ref{thm:vdbound}]\label{thm:Anew}
Let $L=K^1\cup\dots\cup K^m$ be a link with $e$ edges in the $1$-skeleton of a $3$-sphere $S$, and let $e_j$ be the number of edges of $K^j$.
\begin{enumerate}[label=\textup{(\alph*)},nosep]
\item If $S$ is constructible, then $\displaystyle b(L)\le \Big\lfloor\frac{\sum_{j=1}^m e_j}{2}\Big\rfloor=\Big\lfloor\frac{e}{2}\Big\rfloor$.
\item If $S$ is shellable, then $\displaystyle b(L)\le \sum_{j=1}^m\Big\lfloor\frac{e_j}{2}\Big\rfloor$.
\item If $S$ is vertex-decomposable, then $\displaystyle b(L)\le \sum_{j=1}^m\Big\lfloor\frac{e_j}{3}\Big\rfloor$.
\end{enumerate}
\end{thmA}

Of these three results, part (b) is the most surprising and difficult to achieve. Its proof is in the spirit of persistent homology: it records how the link grows alongside a given shelling, and combines a bridge-index estimate of Ehrenborg and Hachimori \cite{EH} with a count of Euler characteristics (Section~\ref{sec:bound}). This argument is specific to shellability, since it is induced by a linear order of the facets. As the next result shows, it does not extend to constructibility.

\begin{thmB}[Theorem~\ref{thm:constructible}]\label{thm:B}
There exists a $3$-dimensional sphere on $13$ vertices that is constructible but unshellable. Its $1$-skeleton contains two disjoint triangles forming a $3$-bridge link whose components are unknotted.
\end{thmB}

By taking suspensions, this example propagates to all dimensions $d \ge 3$ (Corollary~\ref{cor:dims}). Main Theorem~\ref{thm:B} also implies that the number of constructible unshellable $d$-spheres with $N$ facets is exponential in $N$, for every $d\ge3$ (Proposition~\ref{prop:many}).

Finally, to better understand the difference between constructible and shellable spheres, we want to refine Main Theorem~\ref{thm:Anew} even further (and even for knots). This gives us a brand new (un)shellability criterion:

\begin{thmC}[Corollary~\ref{cor:obstruction}]\label{thm:C}
Let $L_1,L_2$ be links in the $1$-skeleton of a $3$-sphere $S$ with $b(L_i)\ge e(L_i)/2$ for $i=1,2$. Suppose some component $K_1$ of $L_1$ and some component $K_2$ of $L_2$ both contain the consecutive pair of edges $(a,b)$ and the consecutive pair $(c,d)$, and that $a$ and $c$ are at even distance along $K_1$ and at odd distance along $K_2$. Then $S$ is not shellable. In particular this holds for two knots $K_1,K_2$ with $b(K_i)\ge e(K_i)/2$.
\end{thmC}

The idea is to look at the case of equality in the count above. If a link is \emph{tight}, i.e.\ it has twice as many edges as its bridge number, then along each component the edges alternate between those that enter the shelling before both of their neighbors on the link and those that enter after both (Theorem~\ref{thm:alt}). If two tight links contain consecutive pairs of edges in positions of opposite parity, the two alternations are incompatible, and the sphere is not shellable.

Main Theorem~\ref{thm:C} has concrete applications. In Appendix~\ref{app:eight} we exhibit five $3$-spheres, $U_{23}$, $U_{24}$, $U_{25}$, $U_{26}$ and $U_{28}$ (the subscript counts the number of vertices), each of which contains two distinct tight knots. All these spheres are proven unshellable by Main Theorem~\ref{thm:C}. We conjecture that none of them is constructible.

\section{Known bounds for knots and links}\label{sec:setting}
Recall that a pure $d$-dimensional simplicial complex is \emph{shellable} if its facets can be ordered $F_1,\dots,F_N$ so that for every $i>1$ the intersection $F_i\cap(F_1\cup\dots\cup F_{i-1})$ is a pure $(d-1)$-dimensional complex; such an order is a \emph{shelling}. It is \emph{constructible} if it is a $d$-simplex, or the union of two constructible $d$-dimensional complexes whose intersection is a constructible pure complex of dimension $d-1$~\cite{HS}. Shellable complexes are constructible.

Let $S$ be a triangulated $3$-sphere with a shelling $F_1,\dots,F_N$, and $B_i=F_1\cup\dots\cup F_i$. For $1\le i<N$, $B_i$ is a combinatorial $3$-ball. For a nonempty face $\sigma$ put $\beta(\sigma)=\min\{i:\sigma\subseteq F_i\}$, so that $\sigma\subseteq B_i$ iff $i\ge\beta(\sigma)$.

A \emph{link} is a disjoint union $L=K^1\cup\dots\cup K^m$ of cycles in the $1$-skeleton of a $3$-sphere $S$, with $e$ edges in total, $o$ of whose components have an odd number of edges; a knot is a link with one component. A component $K$ has vertices $v_1,\dots,v_{e'}$ and edges $g_i=v_iv_{i+1}$ (indices mod $e'$). The \emph{bridge number} $b(L)$ is the minimum, over all positions of $L$ in $S^3$ (that is, embeddings PL isotopic to $L$) relative to a fixed height function $\psi$, of the total number of local maxima of $\psi|_L$; for a knot this is Schubert's bridge number~\cite{Sch}. Call $L$ \emph{tight} if $e=2b(L)$. Since the maxima of $\psi|_L$ are those of its components, $b(L)\ge\sum_i b(K^i)$.

We recall the obstruction to local constructibility of Benedetti and Ziegler~\cite{BZ}, which is in terms of the fundamental group of the complement. A $d$-sphere $S$ is \emph{locally constructible} (LC) if $S-\Delta$ collapses onto a subcomplex of codimension two for some, equivalently every, facet $\Delta$ of $S$~\cite{BZ}; for $d=3$ this means that $S-\Delta$ is collapsible~\cite[Cor.~2.11]{BZ}. Constructible spheres are LC~\cite[Thm.~2.1]{BZ}. The obstruction \cite[Cor.~2.20]{BZ} is stated for arbitrary subcomplexes of codimension two, so it applies verbatim to links:

\begin{proposition}[Benedetti--Ziegler]\label{prop:linkLC}
Let $L$ be a link with $e$ edges in the $1$-skeleton of a $3$-sphere $S$. If $\pi_1(S^3-L)$ has no presentation with $e$ generators, then $S$ is not LC.
\end{proposition}
The rank of $\pi_1(S^3-L)$ is at most the bridge number $b(L)$, since a bridge presentation gives $b(L)$ generators. However, knots with bridge index larger than the rank of their group are known, e.g.\ $8_{16}$, so in general the obstruction of Proposition~\ref{prop:linkLC} is weaker than the following one.

\begin{proposition}[essentially Hachimori--Shimokawa]\label{prop:linkconstr}
Let $L$ be a link with $e$ edges in the $1$-skeleton of a constructible $3$-sphere or $3$-ball. Then $e\ge2b(L)$.
\end{proposition}
\begin{proof}
Hachimori and Shimokawa prove~\cite[Prop.~3.1]{HS}: if $T$ is a tangle (a set of disjoint spanning arcs and knots) in the $1$-skeleton of a constructible regular CW decomposition of a $3$-ball or $3$-sphere, and $i(T)$ is the number of vertices in the relative interior of $T$, then $2b(T)\le2e(T)-i(T)$. In a $3$-sphere a tangle is a link, and for a link $i(L)=e(L)$, since every vertex of a closed component is a relative interior vertex. Hence $2b(L)\le e(L)$.
\end{proof}

For vertex-decomposable complexes there is an even stronger bound:

\begin{proposition}[essentially Ehrenborg--Hachimori]\label{thm:linkvd}
Let $L$ be a link with $e$ edges in the $1$-skeleton of a vertex-decomposable $3$-sphere or $3$-ball. Then $e\ge3b(L)$.
\end{proposition}
For knots this is~\cite[Thm.~5.2]{EH}; for links it is a special case of Theorem~\ref{thm:vdbound}, which we prove in Section~\ref{sec:vd}.

The rest of the paper is dedicated to improving the last two Propositions~\ref{prop:linkconstr} and~\ref{thm:linkvd}. This is achieved in Theorems~\ref{lem:linkbound} and~\ref{thm:vdbound} below.

\section{Proof of Main Theorem A}\label{sec:bound}

In this section, we prove part (b) of Main Theorem~\ref{thm:Anew}. For part (c), see Section \ref{sec:vd}.

Fix a shelling $F_1,\dots,F_N$ of the $3$-sphere $S$, and let $L=K^1\cup\dots\cup K^m$ be a link in its $1$-skeleton, with $e_j$ edges in $K^j$. An edge $g$ of $L$ is \emph{new at step $i$} if $\beta(g)=i$, i.e.\ if $F_i$ is the first facet containing it. Call a component \emph{local} if all its vertices lie in one facet; it is then a cycle of $3$ or $4$ edges contained in the boundary of that facet. Suppose for the moment that no component is local. The new edges at step $i$ then form disjoint paths: a cycle of new edges would be a whole component contained in $F_i$. We call these paths the \emph{new arcs} of step $i$. Each new arc has $k\in\{0,1,2\}$ endpoints on edges of $L$ that were new at earlier steps (\emph{old} endpoints). Since $L$ is a $1$-manifold, the interior vertices of an arc are not old. The last facet $F_N$ has no new edges, because $\partial F_N\subset B_{N-1}$. Let $n_j^k$ be the number of new arcs, over all steps, that lie on $K^j$ and have $k$ old endpoints.

\begin{lemma}\label{lem:count}
Suppose no component of $L$ is local. Then for every $j$ we have:
\begin{enumerate}[label=\textup{(\alph*)},nosep]
\item $n_j^0=n_j^2$;
\item $e_j\ge n_j^0+n_j^1+n_j^2$, with equality if and only if every new arc on $K^j$ is a single edge;
\item $b(L)\le\sum_j n_j^0$.
\end{enumerate}
\end{lemma}
\begin{proof}
(a) Let $K^j_i\subseteq K^j$ be the subgraph formed by the edges new at steps $\le i$ and their endpoints. A new arc with $\ell$ edges and $k$ old endpoints adds $\ell+1-k$ vertices and $\ell$ edges, so it changes the Euler characteristic by $1-k$. Now $K^j_0=\varnothing$ and $K^j_N=K^j$ is a circle, both of Euler characteristic $0$, so $n_j^0-n_j^2=0$.

\noindent (b) The new arcs on $K^j$ have disjoint sets of edges, and every edge of $K^j$ is in exactly one of them. Each arc has at least one edge.

\noindent (c) For $0\le i\le N-1$ let $T_i$ be the union of the edges of $L$ new at steps $\le i$, regarded as a tangle in the ball $B_i$, in the sense of~\cite[Sec.~3]{EH}: a disjoint union of spanning arcs and knots, where a spanning arc has its endpoints on $\partial B_i$ and may run along $\partial B_i$. This is a tangle: a vertex $w$ that ends a path of $T_i$ is followed by an edge of $L$ that is not in $B_i$, so $w$ is not an interior vertex of $B_i$, i.e.\ $w\in\partial B_i$. Also $T_0=\varnothing$, and $T_{N-1}=L$ because $F_N$ has no new edges. We claim that
\begin{equation}\label{eq:step}
b(T_i)\ \le\ b(T_{i-1})+c_i,
\end{equation}
where $a_i,b_i,c_i$ are the numbers of new arcs at step $i$ with $2$, $1$, $0$ old endpoints, respectively, and the bridge index of a tangle is that of~\cite[Def.~3.2]{EH}, and $b(\varnothing)=0$. For $i=1$ the new arcs are disjoint arcs in $\partial F_1$, all with no old endpoint, so they are simultaneously straight and $b(T_1)=c_1$. Let $1<i<N$. Then $B_i=B_{i-1}\cup F_i$ and $B_{i-1}\cap F_i$ is a disk $A_i$, the attaching disk. A new edge of step $i$ is an edge of $\partial F_i$ that is not in $B_{i-1}$, so it meets $B_{i-1}$ only in endpoints; the new arcs are therefore disjoint arcs in the free disk $\overline{\partial F_i\setminus A_i}$, and they touch $A_i$ only in endpoints. An endpoint of a new arc that lies in $B_{i-1}$ either lies on an edge of $T_{i-1}$ (an old endpoint) or is a vertex of $B_{i-1}$ that carries no edge of $T_{i-1}$. In the second case we slide the endpoint of the arc, along $\partial B_i$, off $A_i$ into the free disk; this is an elementary deformation of the tangle in the sense of~\cite[Sec.~3]{EH}, and it does not change $b(T_i)$. After these moves $T_i\cap B_{i-1}=T_{i-1}$, and the closure of $T_i-T_{i-1}$ is a set of $a_i+b_i+c_i$ disjoint arcs in the free disk with $2$, $1$, $0$ points on $T_{i-1}$, respectively. These arcs lie in $\partial F_i\setminus A_i$, hence are simultaneously straight with respect to $A_i=B_{i-1}\cap F_i$, as in the proof of~\cite[Thm.~5.1]{EH}. Now~\cite[Prop.~3.6]{EH} gives~\eqref{eq:step}. Summing~\eqref{eq:step} over $i\le N-1$, and using that the bridge index of the link $T_{N-1}=L$ in the ball $B_{N-1}$ is its bridge number, we get $b(L)\le\sum_ic_i=\sum_jn_j^0$.
\end{proof}

\begin{theorem}\label{lem:linkbound}
A link $L$ with $e$ edges, $o$ of whose components have an odd number of edges, in the $1$-skeleton of a shellable $3$-sphere satisfies $b(L)\le(e-o)/2$. In particular, if a $3$-sphere contains a link with $b(L)>(e-o)/2$, then it is not shellable.
\end{theorem}
\begin{proof}
Let $q$ be the number of local components. A local component bounds a disk $\Delta$ contained in the facet $F$ that contains its vertices, whose interior lies in the interior of $F$: for a triangle this is a face of $F$, and a $4$-cycle in $\partial F$ bounds a disk in $F$ pushed in from the boundary. Since the link lies in the $1$-skeleton, the other components miss this disk, and a small ball around it separates the component from the rest of the link. So it can be added to the others with one local maximum, and $b(L)\le b(L^\circ)+q$, where $L^\circ$ is the union of the non-local components. Each local component has $e_j\in\{3,4\}$ edges, so $\lfloor e_j/2\rfloor\ge1$. By Lemma~\ref{lem:count}, with $n_j^0\le e_j/2$ from (a) and (b), and hence $n_j^0\le\lfloor e_j/2\rfloor$ since $n_j^0$ is an integer,
\begin{equation}\label{eq:chain}
b(L)\ \le\ q+\sum_{K^j\ \text{non-local}}n_j^0\ \le\ q+\sum_{K^j\ \text{non-local}}\Big\lfloor\frac{e_j}2\Big\rfloor\ \le\ \sum_{j}\Big\lfloor\frac{e_j}2\Big\rfloor\ =\ \frac{e-o}2 .\qedhere
\end{equation}
\end{proof}

\medskip

\section{Proof of Main Theorem B}\label{sec:B}

We first build a constructible sphere with $36$ vertices that is not shellable; later we simplify it to a $13$-vertex triangulation.

Let $U$ be the cycle $u_0m_0u_1m_1\cdots u_5m_5u_0$, where the indices are taken modulo six. Call \emph{$i$-th side} the path $u_im_iu_{i+1}$. Write $ij$ for the pair $\{u_i,u_j\}$, and take the three noncrossing perfect matchings
\[
M_X=\{01,25,34\}, \qquad M_Y=\{12,03,45\}, \qquad M_Z=\{23,14,50\}.
\]
For each matching take a disk with boundary $U$. Then draw a chord for each matched pair and subdivide the chord once at a new vertex. Triangulate every one of the four resulting regions by coning its boundary from a new interior vertex. Call the three triangulated disks the \emph{pages} $X,Y,Z$. Since we use disjoint sets of new vertices on different pages, any two pages meet exactly in $U$, and each page induces exactly $U$ on the common vertex set $V(U)$. Inside each page we have two $4$-sided regions, lying between a chord and a side of $U$, and two $8$-sided regions, alternating between two subdivided chords and two subdivided sides. Therefore a page has $12+3+4=19$ vertices and $4+4+8+8=24$ triangles. Call \emph{tips} the centers of the $4$-sided regions, and for $P\in \{X,Y,Z\}$ denote the tips on the even-indexed and odd-indexed sides by $t_P^0$ and $t_P^1$, respectively (Figure~\ref{fig:pages}).

\newcommand{\drawpage}[6]{%
 \begin{scope}[xshift=#1cm]
 \draw[black!60] (0,0) circle (1.12);
 \foreach \i in {0,...,5}{
   \coordinate (u\i) at ({90-60*\i}:1.12);
   \fill[black!55] ({60-60*\i}:1.12) circle (.04);
   \fill (u\i) circle (.045);
   \node[font=\scriptsize,inner sep=0pt] at ({90-60*\i}:1.31) {$\i$};}
 \draw[line width=.8pt] (u#3)--(u#4);
 \fill (0,0) circle (.05);
 \foreach \a in {#5,#6}{
   \pgfmathtruncatemacro{\b}{mod(\a+1,6)}
   \draw[line width=.8pt] (u\a)--({60-60*\a}:.68)--(u\b);
   \fill ({60-60*\a}:.68) circle (.04);}
 \fill ({60-60*#5}:.91) circle (.055);
 \draw[fill=white,line width=.8pt] ({60-60*#6}:.91) circle (.055);
 \node at (0,-1.64) {$#2$};
 \end{scope}}
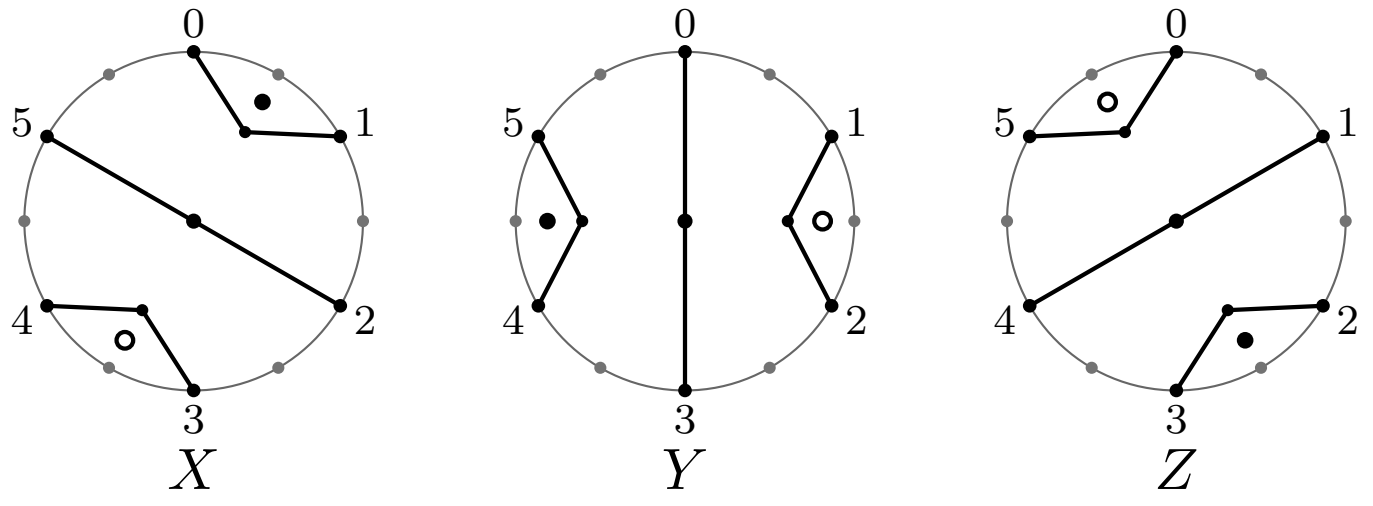
\begin{figure}[htbp]
\centering
\begin{tikzpicture}
\drawpage{0}{X}{2}{5}{0}{3}
\drawpage{3.25}{Y}{0}{3}{4}{1}
\drawpage{6.5}{Z}{1}{4}{2}{5}
\end{tikzpicture}
\caption{The three copies of $U$ are identified by the labels of the vertices. The tips $t_P^0$ and $t_P^1$ are represented filled and empty, respectively. The edges from the region centers to their boundaries are omitted.}
\label{fig:pages}
\end{figure}

Order the pairs of pages cyclically as $(X,Y), (Y,Z), (Z,X)$. For any such pair $(P,Q)$, the pages form a $2$-sphere $P \cup Q$. Moreover, the six chords of $M_P$ and $M_Q$ create a single embedded curve $\gamma_{PQ}$ of twelve edges: indeed, the union of two of the matchings is a disjoint union of cycles of even length, with six vertices in total, and no cycle of length $2$ (a doubled pair), so it is a single $6$-cycle before subdivision.

The curves $U$ and $\gamma_{PQ}$ cross at each $u_i$ because the two chords lie on opposite pages. Hence the sides of $U$ with even index all lie on one side of $\gamma_{PQ}$, and those with odd index on the other. Denote by $H_{PQ}^h$, $h=0,1$, the corresponding closed disks: each contains $t_P^h$ and $t_Q^h$ in its interior, and each contains (on each page) one $4$-sided region and one $8$-sided region. So each $H_{PQ}^h$ has $2(4+8)=24$ triangles, and $H_{PQ}^0\cap H_{PQ}^1=\gamma_{PQ}$. Take a new vertex $g_{PQ}$ and let $G_{PQ}=g_{PQ}*\gamma_{PQ}$ (a disk with twelve triangles). Then, for $h=0,1$, 
\[
J_{PQ}^h=H_{PQ}^h \cup G_{PQ}
\]
is a triangulated $2$-sphere with $36$ triangles, and $J_{PQ}^0\cap J_{PQ}^1=G_{PQ}$. We will fill it by a shellable ball whose interior contains the edge $t_P^ht_Q^h$, except for its endpoints.

\begin{lemma}
\label{lem:cone}
Let $s,t$ be non-adjacent vertices of a triangulated $2$-sphere $J$. Let $D$ be the subcomplex of faces avoiding $s$. Then $s*D$ is a shellable $3$-ball with boundary $J$. Moreover, the edge $st$ is a boundary-parallel spanning arc.
\end{lemma}
\begin{proof}
The closed star of $s$ is a disk bounded by $\lk_J(s)$. Its complementary closed disk is $D$, and $t$ is an interior vertex of $D$ because $s$ and $t$ are non-adjacent. The boundary of $s*D$ is $D\cup (s * \partial D)=J$, and coning a shelling of $D$ gives a shelling of the ball.

Choose a simple path in $D$ from $t$ to a boundary point $z$, all contained in the interior of $D$ except for its endpoint $z$. Then the cone on this path is a disk with boundary $st$ and an arc on $J$, which is the path from $t$ to $z$ followed by $zs$.
\end{proof}

Apply Lemma~\ref{lem:cone} to $J_{PQ}^h$ with $s=t_{P}^h$ and $t=t_{Q}^h$. The tips belong to the interiors of different pages and neither of them lies on $G_{PQ}$, so they are nonadjacent in $J_{PQ}^h$. Write $B_{PQ}^h$ for the resulting cone: exactly four triangles of $J_{PQ}^h$ contain $s=t_{P}^h$, and each of the other $32$ triangles gives one tetrahedron of $B_{PQ}^h$. The choice of the \emph{first} page $P$ of the pair $(P,Q)$ as the cone point is crucial when the sectors are assembled: a face of $B_{PQ}^h$ that is not a face of $J_{PQ}^h$ contains $t_P^h$. Set
\[
W_{PQ}=B_{PQ}^0 \cup B_{PQ}^1, \qquad C_{36}=W_{XY} \cup W_{YZ} \cup W_{ZX},
\]
where all unions are unions of abstract simplicial complexes on the specified vertices.

\begin{lemma}
\label{lem:gluings}
The two halves of $W_{PQ}$ meet exactly in $G_{PQ}$. If $(P,Q)$ and $(Q,R)$ are two of the pairs $(X,Y),(Y,Z),(Z,X)$, then $W_{PQ}\cap W_{QR}$ is exactly the page $Q$. Consequently, $C_{36}$ is a constructible PL $3$-sphere.
\end{lemma}
\begin{proof}
Let $V(\cdot)$ denote vertex sets. Every face of $B_{PQ}^h$ that is not a face of $J_{PQ}^h$ contains the cone point $t_P^h$, and $V(B_{PQ}^h)=V(J_{PQ}^h)$.

The tip $t_P^h$ lies in the interior of $H_{PQ}^h$, so it is not a vertex of $J_{PQ}^{1-h}$. Hence a face common to $B_{PQ}^0$ and $B_{PQ}^1$ contains neither cone point, so it is a face of both $J_{PQ}^0$ and $J_{PQ}^1$, i.e.\ a face of $J_{PQ}^0\cap J_{PQ}^1=G_{PQ}$; conversely $G_{PQ}$ lies in both balls.

We have $V(W_{PQ})=V(P)\cup V(Q)\cup\{g_{PQ}\}$ and $V(W_{QR})=V(Q)\cup V(R)\cup\{g_{QR}\}$. Since $V(P)\cap V(R)=V(U)\subset V(Q)$, a face common to the two sectors has all its vertices in $V(Q)$. Let $\sigma$ be a face of $W_{PQ}$ with $V(\sigma)\subset V(Q)$. The cone points $t_P^h$ of $W_{PQ}$ are not in $V(Q)$, so $\sigma$ is a face of some $J_{PQ}^h$; it does not contain $g_{PQ}$, so it is a face of the sphere $P\cup Q$. A face of the page $P$ with all vertices in $V(U)$ is a face of $U$, because every triangle of a page contains a region center, and an edge of a page joining two vertices of $U$ is a side of $U$. Hence $\sigma$ is a face of the page $Q$. Conversely, the page $Q$ is contained in $H_{PQ}^0\cup H_{PQ}^1$ and in $H_{QR}^0\cup H_{QR}^1$, so it lies in both sectors. This proves that $W_{PQ}\cap W_{QR}$ is the page $Q$.

Each $B_{PQ}^h$ is a shellable $3$-ball with boundary $J_{PQ}^h$ (Lemma~\ref{lem:cone}). The two balls of $W_{PQ}$ meet in the disk $G_{PQ}$, which lies in the boundary of each, so $W_{PQ}$ is a constructible $3$-ball, with boundary $H_{PQ}^0\cup H_{PQ}^1=P\cup Q$. Next, $W_{XY}$ and $W_{YZ}$ meet in the page $Y$, a disk in the boundary of each, so $W_{XY}\cup W_{YZ}$ is a constructible $3$-ball with boundary $X\cup Z$, the union of the pages $X$ and $Z$, which is a $2$-sphere since $X\cap Z=U$. Finally, $W_{ZX}$ is a constructible $3$-ball with boundary $Z\cup X$, and it meets $W_{XY}\cup W_{YZ}$ in $(W_{ZX}\cap W_{XY})\cup(W_{ZX}\cap W_{YZ})=X\cup Z$, which is its whole boundary and a constructible $2$-sphere. So $C_{36}$ is the union of two constructible balls along their whole boundaries, i.e.\ a constructible PL $3$-sphere.
\end{proof}

The vertices of $C_{36}$ are the $12$ vertices of $U$, plus $7$ vertices per page, plus the three vertices $g_{PQ}$, for a total of $36$. The six balls contribute $6\cdot 32=192$ tetrahedra, with no tetrahedron in two of them, since the pairwise intersections are $2$-dimensional. Every triangle belongs to two tetrahedra, so $f_2=384$, and the Euler characteristic of $S^3$ gives $f_1=228$. For $h=0,1$, the three tips of index $h$ span the cycle
\begin{equation}
\label{eq:3-cycle}
K_h=t_X^h\,t_Y^h\,t_Z^h\,t_X^h,
\end{equation}
whose edges $t_X^ht_Y^h$, $t_Y^ht_Z^h$, $t_Z^ht_X^h$ are the edges $st$ in the interiors of $B_{XY}^h$, $B_{YZ}^h$, $B_{ZX}^h$.

Let $L=K_0\cup K_1$. The two triangles $K_0,K_1$ are disjoint, so there is a connected double cover $p\colon M\to|C_{36}|\cong S^3$ branched along $L$, in which every meridian of $L$ interchanges the two sheets. We compute $H_1(M;\mathbb F_3)$ from the decomposition of $C_{36}$ into the three sectors $W_{PQ}$. The preimage of $W_{PQ}$ is denoted $V_{PQ}$, and the preimage of a page $P$ is denoted $\widetilde P$.

\begin{lemma}\label{lem:winding}
Each $V_{PQ}$ is a solid torus. The preimage $\widetilde P$ of the page $P$ is an annulus in $\partial V_{PQ}$, and its core is homologous to $\pm3$ times a generator of $H_1(V_{PQ};\Z)$. The same holds for $Q$.
\end{lemma}
\begin{proof}
The sector $W=W_{PQ}$ is the union of the balls $B^0,B^1$, which meet in the disk $G_{PQ}$, and $L\cap W$ consists of two arcs, $t_P^0t_Q^0\subset B^0$ and $t_P^1t_Q^1\subset B^1$, each of which is a boundary-parallel spanning arc of its ball (Lemma~\ref{lem:cone}). The double cover of a ball branched along an unknotted spanning arc is a ball, and the disk $G_{PQ}$ misses $L$, so it lifts to two disjoint disks $D',D''$.
Gluing the two balls along the first disk joins them into a single ball; gluing along the second disk then adds a $1$-handle. Hence $V_{PQ}$ is a solid torus, and each lift of $G_{PQ}$ is a meridian disk: cutting along one of them leaves a ball.
Its boundary is the double cover of the sphere $S=\partial W=P\cup Q$ branched at the four tips $t_P^0,t_P^1,t_Q^0,t_Q^1$. The page $P$ is a disk containing two branch points, so $\widetilde P$ is an annulus, whose boundary is the preimage of $U=P\cap Q$; it has two components $\widetilde U'$ and $\widetilde U''$, since $U$ separates the branch points two and two.
Let $\mu\in H_1(\partial V_{PQ};\Z)$ be the class of the meridian $\partial D'$, and let $\lambda\in H_1(\partial V_{PQ};\Z)$ be the class of the boundary component $\widetilde U'$ of $\widetilde P$ (equivalently, of the core of $\widetilde P$). The two components $\widetilde U'$ and $\widetilde U''$ are parallel on $\partial V_{PQ}$, so the choice of one of them does not matter. Both $\mu$ and $\lambda$ are primitive classes. We compute their intersection number.

In the four-punctured sphere $P\cup Q\setminus\{\text{tips}\}$, the two curves $U$ and $\gamma=\gamma_{PQ}=\partial G_{PQ}$ meet at the six corners $u_i$. After ignoring the subdivision vertices, their complementary regions are four bigons, one for each $4$-sided region, and four quadrilaterals (the images of the $8$-sided regions after ignoring the subdivision vertices). Each bigon contains one branch point and each quadrilateral contains none. Consequently a bigon lifts to one quadrilateral, and a quadrilateral lifts to two. This gives a quadrangulation of the covering torus with four corners at every vertex. Give each quadrilateral the metric of a unit square and glue by edge isometries. The resulting torus is flat, and the lifts of both curves go straight through each crossing. Thus they realize the geometric intersection of their slopes. There are two parallel lifts in each family, and the six crossings downstairs give twelve upstairs. The four pairs of lifts have the same intersection number, so
\[
4\,|\det(\mu,\lambda)|=12.
\]
This gives the usual slope formula on the four-punctured sphere \cite{Oertel}. In a meridian--longitude basis, the absolute longitude coefficient of $\lambda$ is therefore three, and killing the meridian gives the asserted class in the solid torus.
\end{proof}

\begin{proposition}\label{prop:rank}
The double branched cover $M$ satisfies $\dim_{\mathbb F_3}H_1(M;\mathbb F_3)=2$.
\end{proposition}
\begin{proof}
Take coefficients in $\F_3$ throughout. Write $V_1=V_{XY}$, $V_2=V_{YZ}$, $V_3=V_{ZX}$ and $A=V_1\cup V_2$. The first intersection $V_1\cap V_2=p^{-1}(Y)$ is an annulus. By Lemma~\ref{lem:winding} its core maps to zero in both solid tori. Since the annulus is connected, Mayer--Vietoris gives
\[
H_1(A)\cong H_1(V_1)\oplus H_1(V_2)
\cong\F_3^2 .
\]
The second intersection is all of $T=\partial V_3$, the union of the two annuli above $Z$ and $X$. In the next Mayer--Vietoris sequence consider
\begin{equation}
\label{eq:map}
H_1(T)\longrightarrow H_1(A)\oplus H_1(V_3).
\end{equation}
Choose one of the two lifts of $U$, representing $\lambda$. Its class in $H_1(T)$ is nonzero because it is a primitive slope. It dies in $V_3$ by Lemma~\ref{lem:winding}, and it dies in $A$ since its inclusion factors through $V_1$, where it also represents three times a generator. The map in equation~\eqref{eq:map} has rank at most one. Also, it has rank at least one since its projection to $H_1(V_3)$ is the boundary inclusion of a solid torus and is surjective. Finally, $T$ is connected, so $H_1(M)$ is the cokernel of this rank-one map into $\F_3^3$.
\end{proof}

\begin{theorem}\label{thm:C36}
The complex $C_{36}$ is a constructible $3$-sphere with $f$-vector $(36,228,384,192)$ that is not shellable. The 3-cycles $K_0$ and $K_1$ are vertex-disjoint and unknotted; they form a link $L=K_0\cup K_1$ with $e=6$, $o=2$ and $b(L)=3$. Moreover, the double cover of $S^3$ branched along $L$ has $\dim H_1(\,\cdot\,;\mathbb F_3)=2$.
\end{theorem}
\begin{proof}
By Lemma~\ref{lem:gluings}, $C_{36}$ is a constructible PL $3$-sphere. The cycles $K_0$ and $K_1$ are vertex-disjoint, and neither is a face, because each tetrahedron lies in one ball, and each ball contains tips from only two pages. It remains to show that $b(L)\ge3$. By Proposition~\ref{prop:rank}, $H_1(M;\mathbb F_3)$ has dimension $2$ and so needs two generators. If $L$ has a bridge position with $b$ maxima, then $M$ has a Heegaard splitting of genus $b-1$, so $H_1(M;\mathbb F_3)$ is generated by $b-1$ elements, and $b(L)\ge3$. But since $C_{36}$ is constructible, Proposition~\ref{prop:linkconstr} implies $b(L)\le3$; applied to each component, which has $3$ edges, it gives $b\le1$, so each component is an unknot. Finally, $e=6$ and $o=2$ give $b(L)>(e-o)/2=2$, so $C_{36}$ is not shellable by Theorem~\ref{lem:linkbound}.
\end{proof}

Next, we wish to simplify $C_{36}$ further, to achieve an example with even fewer vertices. It is believed (resp.\ known) that all spheres with $12$ or fewer (resp.\ $10$ or fewer) vertices are shellable (cf.~\cite{Lut}), so $13$ is the smallest number of vertices one can hope for. We achieve it:

\begin{theorem}\label{thm:constructible}
There is a triangulated $3$-sphere $C_{13}$ with $f$-vector $(13,74,122,61)$ that is constructible and not shellable. Its $1$-skeleton contains two disjoint unknotted $3$-cycles $T=\{1,2,3\}$ and $T'=\{4,5,6\}$ that form a link $L$ with $e=6$, $o=2$ and $b(L)=3$.
\end{theorem}

\begin{proof}
Starting from $C_{36}$, we applied random bistellar flips ($2\to3$ and $3\to2$) and edge contractions satisfying the link condition, never removing an edge or a vertex of $L$ and never merging two of its vertices. Contractions alone lead to spheres with $16$ vertices; adding flips with a simulated-annealing acceptance rule, we reached $13$ vertices. A sphere with $14$ vertices, found earlier, is listed in Appendix~\ref{app:C14}. The vertices of $C_{13}$ are $1,\dots,9$ and $a,b,c,d$, for $10,\dots,13$. Its $61$ facets are
\begin{faclist}{14ac 49ac 14bc 14bd 12bd 25bd 479c 47bc 12bc 258d 12cd 1acd 25bc 256c 3acd 28cd 38cd}
14ac 49ac 14bc 14bd 12bd 25bd 479c 47bc 12bc 258d 12cd 1acd 25bc 256c 3acd 28cd 38cd\\
268c 57bc 567c 2567 368c 2678 36ac 1479 1679 1369 1678 1368 1389 369d 389d 1589 1459\\
589d 459a 145a 158a 178a 258a 679c 36ab 35ab 69ab 59ab 59bd 69bd 235a 69ac 46bd 346b\\
347b 346d 357b 347d 2357 237a 37ad 17ad 278a 147d
\end{faclist}
Let $B_1$ be the $3$-ball formed by the first twenty-four tetrahedra, $B_2$ the $3$-ball formed by the next seventeen, and $B_3$ the $3$-ball formed by the last twenty; each of these three blocks of the list is in the order of a shelling of its ball. The complex $B_2\cap B_3$ is a disk, with the $10$ triangles
\par\smallskip\noindent\hfil{\small\ttfamily 69c 69d 59d 59a 36d 25a 28a 78a 17a 147}\hfil\par\smallskip
and $B_1\cap(B_2\cup B_3)$ is a $2$-sphere, with the $22$ triangles
\begin{faclist}{479 49a 79c 9ac 14a 1ad 47b 14d 6ac 67c 36a 3ad 4bd 5bd 57b 368 678}
479 49a 79c 9ac 14a 1ad 47b 14d 6ac 67c 36a 3ad 4bd 5bd 57b 368 678\\
38d 58d 258 257 278
\end{faclist}
Both lists are in the order of a shelling, so $B_2\cup B_3$ is constructible, and then so is $C_{13}=B_1\cup(B_2\cup B_3)$. Being a constructible pseudomanifold without boundary, $C_{13}$ is thus a PL $3$-sphere. Since $e=6$ and $o=2$, by Theorem~\ref{lem:linkbound} it suffices to show $b(L)\ge3$. The moves of the construction above are PL homeomorphisms of the pair (sphere, link), so $(|C_{13}|,L)$ is PL homeomorphic to $(|C_{36}|,K_0\cup K_1)$, and the same holds for the double branched covers. By Proposition~\ref{prop:rank}, the double cover of $S^3$ branched along $L$ has $\dim H_1(\,\cdot\,;\mathbb F_3)=2$, and the argument in the proof of Theorem~\ref{thm:C36} gives $b(L)\ge3$. Finally, $b(L)=3$ and each component is an unknot, by Proposition~\ref{prop:linkconstr} since $C_{13}$ is constructible.
\end{proof}

\begin{remark}\label{rem:92}
It can be verified with SnapPy~\cite{snappy} and Regina~\cite{Reg} that the link $L$ of $C_{13}$ is the link $9^2_{53}$ of Rolfsen's table, up to mirror image.
\end{remark}

$C_{13}$ also solves an open question by Benedetti and Lutz~\cite{BL}:

\begin{proposition}\label{prop:evasive} There exist $3$-balls that are LC but evasive. 
\end{proposition}
\begin{proof}
Let $\Delta'=14bd$ be the fourth facet in the list of facets for $C_{13}$. By the proof of Proposition~\ref{prop:many}, $C_{13}-\Delta'$ is a constructible $3$-ball (with $13$ vertices and $60$ facets). All constructible complexes are LC~\cite[Lemma~2.23]{BZ}. Recall that a complex is non-evasive if it is a point, or if it has a vertex $v$ such that both the link and the deletion of $v$ are non-evasive. We tested this recursively for $C_{13}-\Delta'$, over all choices of the vertex $v$ at each stage, with memoization; the search finds no non-evasive sequence.
\end{proof}

\begin{corollary}\label{cor:dims}
For every $d\ge3$ there is a constructible $d$-sphere that is not shellable.
\end{corollary}
\begin{proof}
By induction on $d$. For $d=3$ this is Theorem~\ref{thm:constructible}. Let $S$ be a constructible $d$-sphere, $d\ge3$, that is not shellable. Suspensions and vertex links preserve both constructibility and shellability. Thus the suspension $\Sigma S=(v*S)\cup(w*S)$ is constructible, and cannot be shellable, since the link of $v$ in $\Sigma S$ is $S$ itself. 
\end{proof}

\begin{remark}\label{rem:sesqui}
\emph{Sesquiconstructible} complexes are constructible complexes in which, in every decomposition $K=U\cup V$ used in the definition, $\partial(U\cap V)$ is empty or again sesquiconstructible. Sesquiconstructibility was introduced by Frohmader~\cite{Fro}, who proved it intermediate between shellability and constructibility. Frohmader also showed that all constructible balls and spheres of dimension at most four are sesquiconstructible~\cite[Prop.~2.8]{Fro}. Hence $C_{13}$ is sesquiconstructible and not shellable. Moreover, the cone over a sesquiconstructible ball or sphere is sesquiconstructible, by induction on the dimension and on the number of facets, using $\partial(v*K)=K\cup v*\partial K$; so the suspension $(v*S)\cup(w*S)$ of a sesquiconstructible sphere $S$, where $\partial S=\varnothing$, is sesquiconstructible. 
Together with Corollary~\ref{cor:dims}, this implies that for every $d\ge3$ there is a sesquiconstructible $d$-sphere that is not shellable. Frohmader proved the Charney--Davis conjecture for odd iterated stellar subdivisions of sesquiconstructible spheres; hence his theorem applies to non-shellable spheres in every dimension $d\ge3$.
\end{remark}

\begin{proposition}\label{prop:many}
Fix $d\ge3$. The number of constructible unshellable $d$-spheres with $N$ facets is exponential in $N$.
\end{proposition}
\begin{proof}
Let $f_d=61\cdot2^{d-3}$. Let $S_d$ be the $(d-3)$-fold suspension of $C_{13}$, that is, $S_3=C_{13}$ and $S_d=(v*S_{d-1})\cup(w*S_{d-1})$ for $d>3$. It is easy to see that $S_d$ is a constructible $d$-sphere with exactly $f_d$ facets. Let $s_d(M)$ be the number of shellable $d$-spheres with $M$ facets, which grows exponentially in $M$~\cite{BZ}. Let $u_d(N)$ be the number of constructible unshellable $d$-spheres with $N$ facets, up to isomorphism. Let $S$ be any shellable $d$-sphere with $M$ facets. Our plan is to construct from $S$ a constructible unshellable $d$-sphere $T_S$ with $M+f_d-2$ facets, and show
\begin{equation}\label{eq:many}
u_d(N)\ \ge\ \frac{s_d(N-f_d+2)}{2\binom N{d+1}}.
\end{equation}

First we claim that there is a facet $\Delta'$ such that $C_{13} - \Delta'$ is constructible. We do this by brute force: with the notation of Theorem~\ref{thm:constructible}, $\Delta'=14bd$ works. In fact, the ball $B_1'=B_1-\Delta'$ is shellable, and so is $B_1'\cap(B_2\cup B_3)$, so $C_{13}-\Delta'$ is a constructible ball with $60$ facets.

Let $\Delta'_d=v*\Delta'_{d-1}$ with $\Delta'_3=\Delta'$, a facet of $S_d$. Then $S_d-\Delta'_d=\bigl(v*(S_{d-1}-\Delta'_{d-1})\bigr)\cup(w*S_{d-1})$ is the union of two cones over constructible complexes, meeting in the constructible complex $S_{d-1}-\Delta'_{d-1}$. By induction, $S_d-\Delta'_d$ is a constructible $d$-ball. Let now $\Delta$ be the last facet in a shelling for the $d$-sphere $S$, so that $S-\Delta$ is a shellable ball. Identify the vertices of $\Delta$ with those of $\Delta'_d$, and give the other vertices of $S$ new names. Put $T_S=(S-\Delta)\cup(S_d-\Delta'_d)$, with $(M-1)+(f_d-1)$ facets. The two balls meet exactly in $\partial\Delta=\partial\Delta'_d$, since a face with all vertices in $\Delta$ is a face of $\Delta$. So $T_S$ is a pseudomanifold without boundary. Since $\partial\Delta$ is a constructible $(d-1)$-sphere, $T_S$ is constructible, hence a sphere. However, no edge of $\Delta'$ is an edge of $T$ or $T'$, so $L$ meets $\Delta'$ only in the vertices $1,4$ and $L$ is still in the $1$-skeleton of $C_{13}-\Delta'$. For $d=3$ this means that $L$ lies in the $1$-skeleton of $T_S$, and $b(L)=3>2=(e-o)/2$, so $T_S$ is not shellable by Theorem~\ref{lem:linkbound}. For $d>3$, let $\sigma$ be the face of $S_d$ formed by the apex $w$ chosen at each suspension step. It is not contained in $\Delta'_d$, whose apexes are the vertices $v$, so its link in $T_S$ is still $C_{13}$. Since links of faces in shellable complexes are shellable, $T_S$ is not shellable.

In $T_S$, the vertex set $\Theta$ of $\Delta$ is one of at most $\binom N{d+1}$ sets, since a $d$-sphere with $N$ facets has at most $N$ vertices. The facets of $T_S$ split into the two strongly connected classes $S-\Delta$ and $S_d-\Delta'_d$, which meet only along $\partial\Theta$, and $S$ is recovered from $T_S$ and $\Theta$ by replacing one of the two classes by the single facet $\Theta$. So $T_S$ determines $S$ up to at most $2\binom N{d+1}$ choices of $(\Theta,\text{class})$, which gives inequality~\eqref{eq:many}. Since $\binom N{d+1}$ is polynomial in $N$, this shows that $u_d(N)$ is at least exponential in $N$.

Since constructible spheres are LC, and there are at most $2^{d^2N}$ LC $d$-spheres with $N$ facets~\cite[Thm.~4.4]{BZ}, we also have $u_d(N)\le 2^{d^2N}$. Hence $u_d(N)$ is exponential in $N$. \end{proof}

\section{Proof of Main Theorem C}\label{sec:shell}

Here we study the tight links, i.e.\ the links with $e=2b(L)$; by Theorem~\ref{lem:linkbound}, in a shellable sphere this is the case of equality with $o=0$. Tight links are especially interesting when the components are linked, so that $b(L)>\sum_i b(K^i)$; in that case the components need not be tight themselves and may even be unknotted.

\begin{definition}\label{def:birth}
Let $g$ be an edge of a link $L$ in the $1$-skeleton of a shelled $3$-sphere, and let $g^-,g^+$ be the two edges of $L$ adjacent to $g$ along its component. We call $g$ a \emph{birth} if $\beta(g)<\beta(g^-)$ and $\beta(g)<\beta(g^+)$, and a \emph{merge} if $\beta(g)>\beta(g^-)$ and $\beta(g)>\beta(g^+)$.
\end{definition}

\begin{theorem}
\label{thm:alt}
Let $L$ be a tight link in the $1$-skeleton of a shellable $3$-sphere, with a fixed shelling. Then every component of $L$ has an even number of edges, and if $g_1,\dots,g_{e'}$ are the edges of a component in cyclic order (indices mod $e'$), then $\beta(g_i)\neq\beta(g_{i+1})$ for all $i$, and the edges alternate between births and merges.
\end{theorem}
\begin{proof}
Let $q$ be the number of local components of $L$. For the non-local components, let $n^k_j$ be as in Lemma \ref{lem:count} and put $n^k=\sum_j n^k_j$, where the sum runs over the non-local components. As in the proof of Theorem~\ref{lem:linkbound}, $b(L)\le q+n^0$, and by Lemma~\ref{lem:count} the $e^\circ$ edges in non-local components satisfy $e^\circ\ge n^0+n^1+n^2=2n^0+n^1$. Since $b(L)=e/2$ and each local component has at least $3$ edges,
\[ \tfrac12\bigl(3q+e^\circ\bigr)\le \tfrac e2=b(L)\le q+n^0\le q+\tfrac12\bigl(e^\circ-n^1\bigr), \]
so $q+n^1\le 0$. Hence there is no local component, $n^1=0$, and all the inequalities above are equalities: $e=2n^0$, and equality in~(b) of the lemma says that every new arc is a single edge.
So every edge $g$ of $L$ is new at step $\beta(g)$, alone, and has either $0$ or $2$ old endpoints. Let $g,g'$ be consecutive edges of a component, sharing the vertex $v$. They are not new at the same step:
if they were, they would form a single new arc with two edges, contradicting the fact that every new arc is a single edge.
If $\beta(g)<\beta(g')$, then $v$ is an old endpoint of $g'$, because $g$ is the other edge of $L$ at $v$; so $g'$ has two old endpoints, and both of its neighbours are new before it, i.e.\ it is a merge. If $g$ has $0$ old endpoints, both of its neighbours are new after it, i.e.\ it is a birth. So an edge is a birth if it has no old endpoint and a merge if it has two; as every edge is one of the two, $\beta$ has a strict local minimum or a strict local maximum at every edge of the cycle, which forces them to alternate, and $e'$ is even.
\end{proof}

\begin{corollary}\label{cor:obstruction}
Let $L_1,L_2$ be links in the $1$-skeleton of a $3$-sphere $S$ with $b(L_i)\ge e(L_i)/2$ for $i=1,2$. Suppose some component $K_1$ of $L_1$ and some component $K_2$ of $L_2$ both contain the consecutive pair of edges $(a,b)$ and the consecutive pair $(c,d)$, and that $a$ and $c$ are at even distance along $K_1$ and at odd distance along $K_2$. Then $S$ is not shellable. In particular, this applies to two knots $K_1,K_2$ with $b(K_i)\ge e(K_i)/2$. (This is Main Theorem~\ref{thm:C} of the introduction.)
\end{corollary}
\begin{proof}
Suppose $S$ is shellable, and fix a shelling. By Theorem~\ref{lem:linkbound}, $b(L_i)\le e(L_i)/2$, so $L_i$ is tight and Theorem~\ref{thm:alt} applies to $K_1$ and $K_2$. Since $(a,b)$ is a pair of consecutive edges of both $K_1$ and $K_2$, we have $\beta(a)\ne\beta(b)$, and $a$ is a birth in $K_i$ if and only if $\beta(a)<\beta(b)$; this does not depend on $i$. The same holds for $c$ and $d$. In $K_1$ the edges $a$ and $c$ are at even distance.
Since the edges of a tight component alternate between births and merges (Theorem \ref{thm:alt}), two edges at even distance have the same type. Hence $\beta(a)<\beta(b)$ if and only if $\beta(c)<\beta(d)$.
In $K_2$ they are at odd distance, so exactly one of them is a birth: $\beta(a)<\beta(b)$ if and only if $\beta(c)>\beta(d)$. These are incompatible, since $\beta(c)\neq\beta(d)$.
\end{proof}

\section{Vertex-decomposability}\label{sec:vd}

For shellable spheres, Theorem~\ref{lem:linkbound} says $b(L)\le(e-o)/2$, where $o$ is the number of components with an odd number of edges; in terms of the numbers $e_j$ of edges of the components, this is $b(L)\le\sum_j\lfloor e_j/2\rfloor$, and $o=\sum_j(e_j\bmod2)$. The same argument that gives $e\ge3b(L)$ in~\cite{EH}, applied to one component at a time, gives the analogous statement for vertex-decomposable spheres, with $2$ replaced by $3$.

\begin{theorem}\label{thm:vdbound}
Let $L=K^1\cup\dots\cup K^m$ be a link with $e$ edges in the $1$-skeleton of a vertex-decomposable $3$-sphere or $3$-ball, let $e_j$ be the number of edges of $K^j$, and let $r=\sum_j(e_j\bmod3)$. Then
\[
b(L)\ \le\ \sum_{j=1}^m\Big\lfloor\frac{e_j}{3}\Big\rfloor\ =\ \frac{e-r}{3}.
\]
\end{theorem}
\begin{proof}
Remove the shedding vertices $x_n,\dots,x_1$ of the vertex-decomposition one at a time, as in~\cite[Prop.~3.4]{EH}. At step $i$ the part of $L$ removed is the set of edges at $x_i$: a straight arc of one or two edges (or nothing), contained in a single component of $L$. Let $a_i,b_i,c_i\in\{0,1\}$ be such that, if the arc is nonempty, exactly one of them equals $1$, namely $a_i$, $b_i$ or $c_i$ according to whether the arc has $2$, $1$ or $0$ endpoints on the part of $L$ that is left; all three are $0$ if the arc is empty. The proof of~\cite[Thm.~5.2]{EH} uses only the estimate~\cite[Prop.~3.6]{EH} for tangles, which holds for links, and gives $b(L)\le\sum_ic_i$. Write $a_i^j,b_i^j,c_i^j$ for the contributions of the arcs lying in $K^j$.

As in Lemma~\ref{lem:count}, an arc with $k$ endpoints on the rest changes the Euler characteristic of the part of $K^j$ removed so far by $1-k$, and this part has Euler characteristic $0$ at the start and at the end, so $\sum_ic_i^j=\sum_ia_i^j$. An arc counted by $a_i^j$ has two edges, and every edge of $K^j$ lies in exactly one arc, so
\[
e_j\ \ge\ \sum_i\bigl(2a_i^j+b_i^j+c_i^j\bigr)\ \ge\ \sum_i\bigl(2a_i^j+c_i^j\bigr)\ =\ 3\sum_ic_i^j .
\]
Hence $\sum_ic_i^j\le\lfloor e_j/3\rfloor$, and summing over $j$ gives $b(L)\le\sum_j\lfloor e_j/3\rfloor$.
\end{proof}

Since $\sum_j\lfloor e_j/3\rfloor\le e/3$, Theorem~\ref{thm:vdbound} implies Proposition~\ref{thm:linkvd}. Let us see an example of a shellable $3$-sphere that Theorem~\ref{thm:vdbound} shows to be not vertex-decomposable, although Proposition~\ref{thm:linkvd} does not.

\begin{proposition}\label{prop:vdexample}
There is a shellable $3$-sphere $S_{16}$ with $16$ vertices and $f$-vector $(16,86,140,70)$ whose $1$-skeleton contains a $4$-cycle and a $5$-cycle forming a link $L$ with $e=9$ and $b(L)=3$, so that $e=3b(L)$. By Theorem~\ref{thm:vdbound}, $S_{16}$ is not vertex-decomposable.
\end{proposition}
\begin{proof}[Proof (computer-assisted)]
Take the sphere $C_{13}$ of Theorem~\ref{thm:constructible}, with vertices written as in its proof, and the link $L_0=T\cup T'$. Subdivide stellarly the edge $12$ of $T$, the edge $45$ of $T'$, and then the edge $56$ of $T'$, with new vertices $e,f,g$ respectively. A stellar subdivision of an edge of $L_0$ is a PL homeomorphism of the pair, so $S_{16}$ is a $3$-sphere containing the same link $L=1\,e\,2\,3\ \cup\ 4\,f\,5\,g\,6$, with $b(L)=3$. The components have $4$ and $5$ edges, so $e=9=3b(L)$, and the bound $e\ge3b$ of Proposition~\ref{thm:linkvd} holds with equality. For vertex-decomposable spheres, however, Theorem~\ref{thm:vdbound} gives $b(L)\le\lfloor4/3\rfloor+\lfloor5/3\rfloor=2$, a contradiction.

It remains to see that $S_{16}$ is shellable. Its $70$ facets, listed here in the order of a shelling found by a SAT solver, are
\begin{faclist}{2bce 1bce 2bde 1bde 2cde 1cde 1acd 17ad 178a 1678 1368 1369 1389 1679 1479 149f 147d}
2bce 1bce 2bde 1bde 2cde 1cde 1acd 17ad 178a 1678 1368 1369 1389 1679 1479 149f 147d\\
14bd 14bc 14ac 14af 49af 49ac 479c 47bc 679c 67cg 57cg 57bc 25bc 25cg 257g 267g 26cg\\
2678 278a 268c 25bd 28cd 258d 258a 158a 1589 159f 15af 69ac 368c 59af 36ac 3acd 38cd\\
389d 369d 37ad 347d 346d 46bd 69bd 69ab 36ab 346b 347b 357b 2357 237a 235a 35ab 59ab\\
59bd 589d
\end{faclist}
\end{proof}

\appendix
\section{Five unshellable spheres with tight knots}\label{app:eight}
Each sphere below satisfies the hypotheses of Corollary~\ref{cor:obstruction} through two knots $K_1,K_2$ with $b(K_i)=e(K_i)/2$, and is therefore not shellable. For $U_{23}$, $U_{24}$, $U_{25}$, $U_{26}$, $U_{28}$ the knots are the hexagons $K_1=1\,2\,3\,4\,5\,6$ and $K_2=1\,2\,3\,6\,5\,4$; the determinants were computed from the facet lists (Fox matrices); all these knots have bridge number $3$ by Schubert's additivity~\cite{Sch}, $b(K\#K')=b(K)+b(K')-1$, since $b(3_1)=b(4_1)=2$; hence $b(K_i)=3=e(K_i)/2$ and all these knots are tight. The common consecutive pairs of edges of $K_1$ and $K_2$ include $(12,23)$ and $(45,56)$, and the edges $12$ and $45$ are at odd distance in $K_1$ and at even distance in $K_2$ (exchange the roles of $K_1$ and $K_2$ in Corollary~\ref{cor:obstruction}). Chirality is not recorded: all identifications are up to mirror image, and ``$\overline{3_1}$'' means a trefoil of the opposite handedness to the other trefoil in the same knot.

\medskip
\noindent{\small\renewcommand{\arraystretch}{1.25}
\begin{tabular}{@{}llllll@{}}
sphere & $f$-vector & $K_1$ & $K_2$ & $e,b$ & $\det K_1,\det K_2$\\\hline
$U_{23}$ & $(23,191,336,168)$ & $3_1\#3_1$ (granny) & $4_1\#4_1$ & $6,3$ & $9,\ 25$\\
$U_{24}$ & $(24,203,358,179)$ & $3_1\#4_1$ & $4_1\#4_1$ & $6,3$ & $15,\ 25$\\
$U_{25}$ & $(25,208,366,183)$ & $4_1\#4_1$ & $4_1\#4_1$ & $6,3$ & $25,\ 25$\\
$U_{26}$ & $(26,221,390,195)$ & $3_1\#3_1$ (granny) & $3_1\#\overline{3_1}$ (square) & $6,3$ & $9,\ 9$\\
$U_{28}$ & $(28,241,426,213)$ & $4_1\#4_1$ & $3_1\#\overline{3_1}$ (square) & $6,3$ & $25,\ 9$\\
\end{tabular}}

\medskip
\noindent Vertices $1,\dots,9$ are written as digits and $10,\dots,35$ as $a,\dots,z$. The knots were identified with SnapPy~\cite{snappy} and knot Floer homology.

\noindent Local constructibility for $3$-spheres is equivalent to admitting $(1,0,0,1)$ as a discrete Morse vector~\cite{BZ}. We verified this with random discrete Morse theory~\cite{BL}: out of $20{,}000$ tries (seed $7$), the optimal vector was reached $33.2\%$, $30.9\%$, $45.7\%$, $28.9\%$ and $31.7\%$ of the times for $U_{23},U_{24},U_{25},U_{26},U_{28}$, respectively. By Corollary~\ref{cor:obstruction} none of the five is shellable. We conjecture that none of them is constructible.

\medskip\noindent The 168 facets of $U_{23}$:
\begin{faclist}{125d 125l 12dl 139j 139k 13ek 13en 13jn 147e 147k 14ek 15cd 15cf 15fi 15hi 15hl 16af}
125d 125l 12dl 139j 139k 13ek 13en 13jn 147e 147k 14ek 15cd 15cf 15fi 15hi 15hl 16af\\
16ai 16fi 178e 178k 189b 189k 18bd 18de 19bj 1afm 1aim 1bdl 1bjn 1bln 1cdn 1cfm 1cim\\
1cin 1den 1hin 1hln 23ek 23en 23kn 248h 248j 249h 249j 256c 256d 25cf 25fi 25il 26ac\\
26af 26di 26fi 28ac 28am 28cg 28eg 28eh 28jm 29cg 29cj 29gh 2afm 2cfm 2cjm 2dil 2egn\\
2ehk 2ghk 2gkn 34ad 34ag 34dg 359g 359n 35ag 35al 35hj 35hl 35jn 368b 368m 36bm 38ac\\
38am 38bd 38cg 38dg 39cg 39cj 39kn 3acl 3adm 3bdm 3chj 3chl 459m 459n 45mn 468b 468j\\
46ad 46aj 46bg 46dg 47ei 47ik 489b 489n 48hn 49bj 49hm 4abg 4abj 4chl 4chm 4cik 4cim\\
4ckl 4eil 4ekl 4hln 4ilm 4lmn 56cd 59gm 5abg 5abj 5ajk 5akl 5bgm 5bjn 5bmn 5ehi 5ehk\\
5eil 5ekl 5hjk 68jm 6ack 6adi 6ajk 6bgm 6cdk 6dgk 6gjk 6gjm 78ei 78ik 89kn 8deg 8ehi\\
8hin 8ikn 9ghm ackl adim bdlm blmn cdkn chjm cikn degn dgkn dilm ghjk ghjm
\end{faclist}

\medskip\noindent The 179 facets of $U_{24}$:
\begin{faclist}{1279 127a 129a 137a 137d 13ak 13dl 13jk 13jo 13lo 14dh 14dl 14hl 15ah 15ak 15bn 15bo}
1279 127a 129a 137a 137d 13ak 13dl 13jk 13jo 13lo 14dh 14dl 14hl 15ah 15ak 15bn 15bo\\
15ho 15jk 15jn 16bj 16bn 16jn 179h 17dh 19ah 1bjo 1hlo 23jl 23jo 23lo 24dh 24dk 24hl\\
24jl 24jo 24ko 25cf 25ci 25fg 25gi 269a 269e 26ab 26bn 26eg 26gn 279b 27ab 28ci 28cm\\
28df 28dh 28ef 28eg 28gi 28hm 29be 2ben 2cdf 2cdm 2dkm 2efn 2fgn 2hlo 2hmo 2kmo 34ce\\
34ci 34ei 3578 357g 358m 35cg 35ci 35im 368f 368h 36fh 378a 37dg 38ag 38ef 38eg 38hm\\
3adg 3adk 3ceg 3dkl 3efh 3ehi 3him 3jkl 4578 457m 458m 478a 47ab 47bk 47fk 47fm 48ai\\
48ci 48cm 4abj 4aei 4aem 4afj 4afm 4bjo 4bko 4cem 4dkl 4fjk 4jkl 569e 569o 56eo 57gi\\
57im 59be 59bo 5acf 5ack 5afj 5ahn 5ajn 5ben 5cgk 5ehn 5eho 5fgk 5fjk 67dh 67do 67fh\\
67fo 68df 68dh 69an 69no 6abj 6ajn 6cde 6cdf 6ceg 6cfo 6cgn 6cno 6deo 79bk 79hk 7dgi\\
7dio 7fhk 7fmo 7imo 8agi 9ahn 9bko 9hkn 9kno acfo acko adei adem adgi adkm afmo akmo\\
cdem cgkn ckno deio efhn ehio fgkn fhkn himo
\end{faclist}

\medskip\noindent The 183 facets of $U_{25}$:
\begin{faclist}{1257 125p 127p 138d 138n 13dm 13fi 13fn 13io 13mo 14fi 14fn 14in 157d 15de 15em 15mp}
1257 125p 127p 138d 138n 13dm 13fi 13fn 13io 13mo 14fi 14fn 14in 157d 15de 15em 15mp\\
16bd 16be 16de 17dp 18do 18no 1bdm 1bem 1dop 1ino 1mop 23im 23io 23mo 246c 246e 247n\\
247p 24cn 24ep 257k 25de 25dk 25ep 26bd 26bh 26ch 26de 27kn 2bcl 2bcn 2bdn 2bhj 2bjo\\
2blo 2chl 2dkn 2hjl 2ijl 2ijo 2ilm 2lmo 349a 349g 34ag 357h 357j 35ag 35ah 35bj 35bl\\
35gl 36jk 36jp 36kp 378c 378j 37ch 38cd 38jn 39al 39gl 3ahl 3bcl 3bcp 3bjk 3bkp 3cdj\\
3chl 3cjp 3djn 3dmn 3fim 3fmn 458l 458o 45lo 467h 467n 46be 46bh 46cm 46im 46in 47dh\\
47dp 48dg 48do 48gp 48lp 49ad 49dk 49gk 4adg 4bep 4bhk 4bkp 4cmn 4dhk 4dop 4fim 4fmn\\
4gkp 4lop 56ci 56cm 56im 57aj 57ak 57dh 58ln 58no 5agj 5ahk 5bjo 5blo 5cem 5cep 5cip\\
5dhk 5gjo 5gln 5gno 5imp 678c 678l 67ch 67kl 67kn 68cg 68gp 68lp 6cgi 6gin 6gkn 6gkp\\
6jkl 6jlp 78aj 78al 7akl 8ajn 8aln 8cdg 9adn 9aln 9dkn 9gkn 9gln adgj adjn ahkl bcem\\
bcep bcmn bdmn bhjk cdgj cgij cijp gijo gino hjkl ijlp ilmp lmop
\end{faclist}

\medskip\noindent The 195 facets of $U_{26}$:
\begin{faclist}{125c 125q 12cq 13in 13iq 13km 13kq 13mo 13no 14jo 14jp 14op 15af 15ag 15cg 15fq 168a}
125c 125q 12cq 13in 13iq 13km 13kq 13mo 13no 14jo 14jp 14op 15af 15ag 15cg 15fq 168a\\
168g 16ag 18af 18cg 18ci 18fi 19km 19kn 19mo 19np 19op 1ciq 1fil 1flq 1iln 1jno 1jnp\\
1kln 1klq 23ai 23an 23in 246b 246g 24bp 24gp 256b 256c 25bg 25gq 26cq 26gq 2adi 2adl\\
2aln 2bgm 2bhl 2bhp 2blm 2dhi 2dhm 2dlm 2gmp 2hil 2hmp 2iln 3489 348d 349d 3578 357k\\
3589 359d 35de 35em 35ko 35mo 36hj 36hp 36jp 378c 37cm 37km 38cd 3aiq 3anq 3bgj 3bgm\\
3bhl 3bhp 3bjp 3blm 3cde 3cel 3clm 3egl 3egm 3ghj 3ghl 3koq 3noq 457i 457k 45ik 469e\\
469k 46bf 46ce 46ci 46fk 46gi 47ci 47cm 47km 489l 48dl 49dm 49el 49km 4bfj 4bjp 4cel\\
4clm 4dlm 4fik 4fio 4fjo 4gio 4gop 56bc 578i 589o 58ik 58ko 59dm 59mo 5afn 5agq 5anq\\
5bcg 5def 5dfh 5dhm 5efn 5ehm 5ehn 5fhq 5hnq 68ad 68dg 69ep 69kn 69np 6adi 6agq 6aiq\\
6bcf 6cef 6ciq 6dgi 6efn 6ehn 6ehp 6fkn 6hjn 6jnp 78ci 89lo 8adl 8afk 8akl 8cdg 8fik\\
8klo 9elp 9lop afkn akln bcfg bfgj cdef cdfg dfgh dghi eglp egmp ehmp fghj fhjq filo\\
fjoq floq ghil gilo glop hjnq jnoq kloq
\end{faclist}

\medskip\noindent The 213 facets of $U_{28}$:
\begin{faclist}{12ag 12ao 12go 13dp 13dq 13hq 13hr 13ik 13ir 13kp 14ck 14cp 14kp 157l 157q 15hl 15hr}
12ag 12ao 12go 13dp 13dq 13hq 13hr 13ik 13ir 13kp 14ck 14cp 14kp 157l 157q 15hl 15hr\\
15jq 15jr 16jq 16jr 16qr 17ah 17ao 17hq 17lo 1afg 1afh 1ckn 1cnp 1dnp 1dnq 1fgl 1fhl\\
1glo 1ikq 1iqr 1knq 23ik 23ir 23kr 24dh 24dl 24hn 24lq 24nq 256e 256g 257e 257j 25gj\\
267j 267m 26ag 26ao 26em 26jo 278b 278e 27bd 27dm 28bm 28em 2bdl 2blm 2dhm 2gjo 2hlm\\
2hln 2ikq 2iqr 2knq 2knr 2lnr 2lqr 349m 349o 34mo 359f 359o 35cf 35cs 35mo 35ms 368b\\
368n 36bn 378b 378n 37bd 37dm 37gm 37gn 39fm 3acf 3acq 3adf 3adp 3apq 3bdg 3bgn 3cdq\\
3cds 3dfm 3dgs 3gms 3hkp 3hkr 3hpq 458a 458d 45ad 467a 467c 468a 468n 46ch 46hn 47ah\\
47cp 47hp 48dn 49do 49dq 49mq 4adh 4chk 4dlo 4dnq 4hkp 4loq 4moq 56eg 57eq 57jr 57lr\\
58am 58dp 58ms 58ps 59bc 59bf 59cs 59os 5adp 5amp 5bch 5bfn 5bhn 5cfn 5chk 5ckn 5egj\\
5ejq 5hkr 5hln 5knr 5lnr 5mop 5ops 67ao 67cg 67gm 67jr 67or 68ac 68bc 6acg 6bch 6bhn\\
6egm 6joq 6oqr 78ep 78np 7cgn 7cnp 7epq 7hpq 7lor 8acm 8bcm 8dnp 8egm 8egp 8gms 8gps\\
9bcm 9bfm 9cdq 9cds 9cmq 9dos acfg acmq adfh ampq bdgl bfgl bfgn bflm cfgn dfhm dglo\\
dgos egjp ejpq fhlm gjop gops jopq loqr mopq
\end{faclist}

\section{Another unshellable constructible sphere}\label{app:C14}
Before $C_{13}$ we found a constructible unshellable sphere $C_{14}$ with $14$ vertices and $f$-vector $(14,77,126,63)$. Its vertices are $1,\dots,9$ and $a,\dots,e$ for $10,\dots,14$, and its $63$ facets are

\medskip
\begin{faclist}{257a 157a 27ab 279b 29be 15ac 24ab 16ac 17ab 39be 179b 139b 16ab 1389 46ab 138b 168b}
257a 157a 27ab 279b 29be 15ac 24ab 16ac 17ab 39be 179b 139b 16ab 1389 46ab 138b 168b\\
468b 468a 1689 348b 689a 348a 69ac 389a 59ac 359a 359e 2bce 12ce 124e 249e 124c 24bc\\
249d 145e 145c 459e 459c 157e 16ce 34bc 3bce 36ce 367e 36cd 167e 1679 34cd 357e 34ad\\
2357 2367 2679 2356 356d 256d 35ad 25ad 24ad 269d 69cd 49cd
\end{faclist}

\medskip\noindent As for $C_{13}$, the first $28$, the next $13$ and the last $22$ facets form shellable balls $B_1,B_2,B_3$, with $B_2\cap B_3$ a shellable disk and $B_1\cap(B_2\cup B_3)$ a shellable $2$-sphere, so $C_{14}$ is constructible. The triangles $T=\{1,2,3\}$ and $T'=\{4,5,6\}$ form a link $L$ with $e=6$, $o=2$ and $b(L)=3$, since the double cover branched along $L$ has $\dim H_1(\,\cdot\,;\mathbb F_3)=2$; hence $C_{14}$ is not shellable by Theorem~\ref{lem:linkbound}. In fact, $L$ is the pretzel link $P(3,3,3,-1)$, up to mirror image.

\subsection*{Acknowledgements and AI disclosure}
Claude Sonnet 5.5 was used in the preparation of some parts of this work, to identify the knots or links, and to perform the computations to obtain the example $C_{13}$ from $C_{36}$ (which was human generated). The examples in the appendix were also generated by Claude Sonnet 5.5. All the material was reviewed and verified by the authors, who take full responsibility for the proofs and the main results.
BB is supported by Simons gift MPS-TSM-00002873; GM is supported by Horizon Europe ERC Grant number: 101045750 / Project acronym: HodgeGeoComb.

\end{document}